\documentclass[11pt]{article} 
\begin{document} 
\newtheorem{prop}{Proposition}[section]
\newtheorem{Def}{Definition}[section] \newtheorem{theorem}{Theorem}[section]
\newtheorem{lemma}{Lemma}[section] \newtheorem{Cor}{Corollary}[section]

\title{\bf An improved local well-posedness result for the one-dimensional 
Zakharov system}
\author{{\bf
 Hartmut Pecher}\\
 Fachbereich Mathematik und Naturwissenschaften\\
  Bergische Universit\"at Wuppertal\\
 Gau{\ss}str.  20
  \\ D-42097 Wuppertal\\
   Germany\\
    e-mail Hartmut.Pecher@math.uni-wuppertal.de}
     \date{}
 \maketitle

\begin{abstract}
The 1D Cauchy problem for the Zakharov system is shown to be locally well-posed 
for low regularity Schr\"odinger data $u_0 \in \widehat{H^{k,p}}$ and wave data 
$(n_0,n_1) \in \widehat{H^{l,p}} \times \widehat{H^{l-1,p}}$ under certain 
assumptions on the parameters $k,l$ and $1<p\le 2$, where 
$\|u_0\|_{\widehat{H^{k,p}}} := \| \langle \xi \rangle^k 
\widehat{u_0}\|_{L^{p'}}$ , generalizing the results for $p=2$ by Ginibre, 
Tsutsumi, and Velo. Especially we are able to improve the results from the 
scaling point of view, and also allow suitable $k<0$ , $l<-1/2$ , i.e. data 
$u_0 \not\in L^2$ and $(n_0,n_1)\not\in H^{-1/2}\times H^{-3/2}$, which was 
excluded in the case $p=2$.
\end{abstract}
\renewcommand{\thefootnote}{\fnsymbol{footnote}}
\footnotetext{\hspace{-1.8em}{\it 2000 Mathematics Subject Classification:} 
35Q55, 35L70 }
\normalsize
 
\setcounter{section}{0}
\section{Introduction and main results} Consider the (1+1)-dimensional Cauchy 
problem for the 
Zakharov system
\begin{eqnarray}
\label{1}
 iu_t + u_{xx} & = & nu
\\ \label{2}
n_{tt}-n_{xx} & = & (|u|^2)_{xx}
\\
\label{3}
 u(0) \quad = \quad u_0 \quad , \quad
n(0) \,& = &\, n_0 \quad , \quad n_t(0) \quad = \quad n_1
 \end{eqnarray}
where $u$ is a complex-valued und $n$ a
real-valued function defined for $(x,t) \in {\bf R}\times{\bf R}^+$.

The Zakharov system was introduced in \cite{Z} to describe Langmuir turbulence 
in a plasma.

The Zakharov system (\ref{1}),(\ref{2}),(\ref{3}) can be 
transformed into a first order system in $t$ as follows: With
$$ n_{\pm} := n \pm iA^{-1/2} n_t \, , {\mbox i.e.} \, n = \frac{1}{2}(n_++n_-) 
\, , \, 2iA^{-1/2}n_t = n_+ - n_- \, , \, A:= -\partial_x^2 $$
we get
\begin{eqnarray}
\label{3.1}
iu_t + u_{xx} & = & \frac{1}{2}(n_+ + n_-)u \\
\label{3.2}
in_{\pm t} \mp A^{1/2} n_{\pm} & = & \pm A^{1/2}(|u|^2) \\
\label{3.3}
u(0) = u_0 \quad , \quad n_{\pm}(0) & = & n_0 \pm iA^{-1/2} n_1 =: n_{\pm 0} 
\, .
\end{eqnarray}
This problem was considered for data in $L^2$-based Sobolev spaces in detail in 
the last decade, especially low regularity local well-posedness results were 
given by Ginibre, Tsutsumi, and Velo \cite{GTV} for data $u_0 \in H^k $ , $ n_0 
\in H^l$ , $n_1 \in H^{l-1}$ under the following assumptions:
$$ -\frac{1}{2} < k-l \le 1 \quad , \quad 2k \ge l + \frac{1}{2} \ge 0 \, . $$
In this paper there were also given corresponding results in arbitrary space 
dimension. It was also shown that these results are sharp within the used 
method, namely the Fourier restriction norm method initiated by Bourgain and 
Klainerman-Machedon and further developed by Kenig, Ponce, Vega and others. It 
could 
also be shown by Colliander, Holmer, and Tzirakis \cite{CHT}, that global 
well-posedness in the case $k=0$ , $ l=-1/2$ holds true. Holmer \cite{H} was 
able to show that the one-dimensional local well-posedness theory is sharp in 
the sense that the problem is locally ill-posed in some cases, where the 
assumptions on $k,l$ in \cite{GTV} are violated, more precisely: if $0<k<1$ and 
$2k>l+1/2$ , or, if $k \le 0$ and $l>-1/2$, or, if $k=0$ and $l<-3/2$. 
Moreover, 
the mapping data upon solution is not $C^2$, if $k \in {\bf R}$ , $ l <-1/2 $. 
Ill-posedness for $k<0$ and $l\le -3/2$ was shown by Biagioni and Linares 
\cite{BL}.

The minimal values $k=0$ , $l=-1/2$ are far from critical, if one compares them 
with those being critical for a scaling argument, namely $k=-1$ and 
$l=-3/2$ . The heuristic scaling argument here is the following (for details we 
refer to 
\cite{GTV}): Ignoring the term $A^{\frac{1}{2}}n_{\pm}$ in equation (\ref{3.2}) 
the system (\ref{3.1}),(\ref{3.2}),(\ref{3.3}) is invariant under the dilation
\begin{eqnarray}
\label{aa}
u(x,t) & \longrightarrow & u_{\mu}(x,t) = \mu^{\frac{3}{2}} u(\mu x,\mu ^2 t) \\
\label{bb}
n_{\pm}(x,t) & \longrightarrow & n_{\pm \mu}(x,t) = \mu ^2 n_{\pm}(\mu x, \mu ^2 
t) \, .
\end{eqnarray}
Because
\begin{equation}
\label{a}
\|u_{\mu}(x,0)\|_{\dot{H}^k} = \mu ^{k+1} \|u_0\|_{\dot{H}^k}
\end{equation}
and
\begin{equation}
\label{b}
\|n_{\pm \mu}(x,0)\|_{\dot{H}^l} = \mu ^{l+\frac{3}{2}} \|n_{\pm 
0}\|_{\dot{H}^l}
\end{equation}
the system is critical for $k=-1$ and $l=-\frac{3}{2}$ . If namely the lifespan 
of $(u,n_+,n_-)$ were $T$ the lifespan of $(u_{\mu},n_{+ \mu},n_{- \mu})$ would 
be $T \mu^{-2}$. 
So, if $k<-1$ or $l<-\frac{3}{2}$, one would have both the norm of the data and 
the lifespan of the solution $(u,n_+,n_-)$ going to zero as $\mu \to \infty$ , 
which strongly indicates ill-posedness.

It is interesting to compare the situation with the corresponding problem for 
the cubic Schr\"odinger equation
\begin{equation}
\label{NLS}
iu_t + u_{xx} + |u|^2u = 0 \quad , \quad u(0) = u_0 \, 
\end{equation}
which is known to be (globally) well-posed for data $u_0 \in H^s$ , $s \ge 0$ 
\cite{Y} (cf. also \cite 
{CW}) , and locally ill-posed for $s<0$ \cite{KPV3}, whereas scaling 
considerations suggest as the critical value $s=-1/2$. This problem is of 
special interest also for the Zakharov system, because the cubic Schr\"odinger 
equation is the formal limit for $c \to \infty$ of the Zakharov system modified 
by replacing $\partial_t^2 - \partial_x^2$ by $c^{-2} \partial_t^2 - 
\partial_x^2$.
Now, for nonlinear Schr\"odinger equations it was suggested to leave the 
$H^s$-scale of the data by Cazenave, Vega, and Vilela \cite{CVV} and Vargas and 
Vega \cite{VV}. For the cubic Schr\"odinger equation local (and even global) 
well-posedness has been shown for data with infinite $L^2$-norm. A. Gr\"unrock 
\cite{G2} was able to show in this case local well-posedness for data $u_0 \in 
\widehat{H^{s,r}}$ , where
$$ \|u_0\|_{\widehat{H^{s,r}}} := \| \langle \xi \rangle^s \widehat{u_0} 
\|_{L_{\xi}^{r'}} \quad , \quad 1/r + 1/r' = 1 \, , $$
if $ s \ge 0$ and $1<r<\infty$ . Moreover, he could show global well-posedness 
for $2\ge r \ge 5/3 $ , $u_0 \in \widehat{H^{0,r}}$ , and also local 
ill-posedness for the cubic Schr\"odinger equation in 
$\widehat{H^{s,r}}$ for any $1<r<\infty$ and $ -1/r' < s < 0$ . The 
well-posedness results were proven by a modified Fourier restriction norm 
method 
(for $p \neq 2$), which was developed by A. Gr\"unrock in \cite{G1}, where 
these 
ideas were applied to the modified KdV equation.

The aim of the present paper is to prove local well-posedness results for the 
Zakharov system with data $u_0 \in \widehat{H^{k,p}}$ , $ n_0 \in 
\widehat{H^{l,p}}$ , $ n_1 \in \widehat{H^{l-1,p}}$ under suitable assumptions 
on $k,l,p$ , which allow to weaken the assumptions on the data from the scaling 
point of view, thus improving the $L^2$-based results in this sense, and also 
allow to get results for certain data $u_0 \not\in L^2$ 
and $(n_0,n_1) \not\in H^{-1/2} \times H^{-3/2}.$ Details are given in section 
2 and 3. Especially we can show that local well-posedness holds for data 
$(u_0,n_0,n_1) \in \widehat{H^{k,p}} \times \widehat{H^{l,p}} \times 
\widehat{H^{l-1,p}}$ for suitable $k<0,$ $l<-1/2$ , and $1<p<2$ in contrast to 
the above-mentioned ill-posedness results of Holmer \cite{H} for the Zakharov 
system, and also in contrast to Gr\"unrock's ill-posedness results for the 
cubic 
Schr\"odinger equation, so that the limit of the c-dependent Zakharov system as 
$ c \to \infty$ must be singular. We are also able to choose $k=0$ and $l>-1/2$ 
, a choice which was not possible in the $L^2$-case (cf. \cite{H} again).

We prove our results by a modification of the  Fourier restriction norm method,
originally due to J. Bourgain \cite{B1},\cite{B2}, and derive the crucial 
estimates for the nonlinearities using a variant of the Schwarz method 
introduced by Kenig, Ponce and Vega \cite{KPV1},\cite{KPV2} adapted to the 
$L^p$-theory. In principle these estimates are proven along the lines of 
\cite{GTV}. 

We recall the modified Fourier restriction norm method in the following. For 
details we refer to the paper of A. Gr\"unrock (cf. \cite{G1}, Chapter 2).
Our solution spaces are the Banach spaces
$$ X_r^{l,b} := \{ f \in {\cal S}'({\bf R}^2) : \|f\|_{X^{l,b}_r} < \infty\} \, 
, $$
where $l,b \in {\bf R}$ , $1<r<\infty$ , $1/r + 1/r' = 1$ and
$$ \|f\|_{X^{l,b}_r} := \left( \int d\xi d\tau \langle \xi \rangle^{lr'} 
\langle 
\tau + \phi(\xi) \rangle^{br'} |\hat{f}(\xi,\tau)|^{r'} \right)^{1/r'} \, , $$
where $\phi: {\bf R} \to {\bf R} $ is a given smooth function of polynomial 
growth. The dual space of $X^{l,b}_r$ is $X^{-l,-b}_{r'}$ , and the Schwartz 
space 
is dense in $X^{l,r}_p$ . The embedding $X^{l,b}_r \subset C^0({\bf 
R},\widehat{H^{l,r}})$ is true for $b>1/r$ . We have
$$ \|f\|_{X^{l,b}_r} = \left(\int d\xi d\tau \langle \xi \rangle^{lr'} \langle 
\tau \rangle^{br'} \left| {\cal F}(e^{-it\phi(-i\partial_x)}f)(\xi,\tau) 
\right|^{r'} \right)^{1/r'} $$
and
$$ \| \psi e^{it\phi(-i\partial_x)} u_0\|_{X^{l,b}_r} \le c_{\psi} 
\|u_0\|_{\widehat{H^{l,r}}} $$
for any $\psi \in C^{\infty}_0({\bf R}_t)$ .

If $v$ is a solution of the inhomogeneous problem
$$ iv_t - \phi(-i\partial_x)v = F \, , \, v(0) = 0 $$
and $\psi \in C^{\infty}_0({\bf R}_t)$ with $supp \, \psi \subset (-2,2) $ , $ 
\psi \equiv 1 $ on $[-1,1]$ , $ \psi(t) = \psi(-t) $ , $ \psi(t) \ge 0 $ , $ 
\psi_{\delta}(t) := \psi(\frac{t}{\delta}) $ , $ 0<\delta \le 1$ , we have for 
$1<r<\infty$ , $ b'+1 \ge b \ge 0 \ge b' > - 1/r' $ :
$$ \|\psi_{\delta} v\|_{X^{l,b}_r} \le c \delta^{1+b'-b} \|F\|_{X^{l,b'}_r} \, 
. 
$$
For the reduced wave part $\phi(\xi) = \pm |\xi|$ we use the notation 
$X^{l,b}_{\pm,r}$ instead of $X^{l,b}_r$ , whereas for the Schr\"odinger part 
$\phi(\xi) = \xi^2$ we simply use $X^{l,b}_r$ . We also use the localized 
spaces
$$ X^{l,b}_r(0,T) := \{ f = \tilde{f}_{|[0,T] \times {\bf R}} : \tilde{f} \in 
X^{l,b}_r \} \, , $$
where
$$ \|f\|_{X^{l,b}_r(0,T)} := \inf \{ \|\tilde{f}\|_{X^{l,b}_r} : f 
= \tilde{f}_{|[0,T]\times{\bf R}}\} \, .$$
Especially we use \cite{G1}, Theorem 2.3, which we repeat for convenience. 
\begin{theorem}
\label{Theorem 0.1}
Consider the Cauchy problem
\begin{equation}
\label{0*}
u_t - i \phi(-i\partial_x) u = N(u) \quad , \quad u(0) = u_0 \in 
\widehat{H^{s,r}} \, , 
\end{equation}
where $N$ is a nonlinear function of $u$ and its spatial derivatives. Assume 
for 
given $s \in {\bf R} $ , $1<r<\infty$ , $\alpha \ge 1$ there exist $ b > 1/r $ 
, 
$b-1 < b' \le 0$ such that the estimates
$$ \|N(u)\|_{X^{s,b'}_r} \le c \|u\|_{X^{s,b}_r}^{\alpha} $$
and
$$ \|N(u)-N(v)\|_{X^{s,b'}_r} \le c (\|u\|_{X^{s,b}_r}^{\alpha -1} + 
\|v\|_{X^{s,b}_r}^{\alpha -1}) \|u-v\|_{X^{s,b}_r}$$
are valid. Then there exist $T = T(\|u_0\|_{\widehat{H^{s,r}}}) > 0 $ and a 
unique solution $u \in X^{s,b}_r[0,T]$ of (\ref{0*}). This solution belongs to 
$C^0([0,T],\widehat{H^{s,r}})$ , and the mapping $u_0 \mapsto u$ , 
$\widehat{H^{s,r}} \to X^{s,b}_r(0,T_0)$ is locally Lipschitz continuous for 
any 
$T_0 < T$ .
\end{theorem}

The main result of this paper is the following
\begin{theorem}
\label{Theorem}
Let $1<p\le 2$ , $ \frac{1}{p} + \frac{1}{p'} = 1$ , $1 \ge b,b_1 > \frac{1}{p} 
$ .
\begin{itemize}
\item In the case $k\ge 0$ assume
$$ l \ge -\frac{1}{p} \, , \, k-l < 2(1-b_1) \, , \, l \le 2k-\frac{1}{p'} \, , 
\, l+1-k < \frac{1}{p} +2(1-b) \, , \, l+1-k \le 2b_1 \, .$$
\item In the case $k<0$ assume
$$ k \ge - \frac{1}{p} \, , \, l \ge -\frac{1}{p} \, , \,
l+k > \frac{1}{p} - 2b_1 \, , l+k > \frac{1}{p} - 2b \, , \, l+k > -\frac{1}{p} 
-2(1-b_1) \, , $$
 $$ k-l < 2(1-b_1) \, , \, 2k > \frac{1}{p}-b_1 \, , \, 2k \ge l+\frac{1}{p'} 
\, , \, 2k > -(1-b) \, . $$
\end{itemize}
Let $u_0 \in \widehat{H^{k,p}}$ , $ n_{\pm 0} \in \widehat{H^{l,p}} $ . Then 
the Cauchy problem (\ref{3.1}),(\ref{3.2}),(\ref{3.3}) is locally well-posed, 
i.e. there exists a unique local solution $u \in X^{k,b_1}_p(0,T)$ , $ n_{\pm} 
\in X^{l,b}_{\pm,p}(0,T)$. This solution satisfies $u \in 
C^0([0,T],\widehat{H^{k,p}})$ , $ n_{\pm} \in C^0([0,T],\widehat{H^{l,p}})$ , 
and the mapping data upon solution is locally Lipschitz continuous.
\end{theorem}

The estimates for the nonlinearities are given in section 3 and the short proof 
of this theorem as a consequence of these estimates in section 4.

{\bf Remark:} The assumption $n_{\pm 0} \in \widehat{H^{l,p}}$ requires 
$n_0,A^{-1/2}n_1 \in \widehat{H^{l,p}}$. This last assumption on $n_1$ can also 
be replaced by the condition $n_1 \in \widehat{H^{l-1,p}}$.  One way to see 
this is to modify the transformation of the original Zakharov system into the 
first order system in $t$ as follows: replace the wave equation by 
$n_{tt}-n_{xx}+n=(|u|^2)_{xx} +n$ and define $n_{\pm}:=n \pm i \tilde{A}^{-1/2} 
n_t$, where $\tilde{A} := -\partial_x^2 +1.$ This leads to the modified reduced 
wave equation:
$$ in_{\pm t} \mp \tilde{A}^{1/2}n_{\pm} = \pm A \tilde{A}^{-1/2}(|u|^2) \mp 
(1/2) \tilde{A}^{-1/2}(n_+ + n_-) \, . $$
Now it is easy to see that this modified nonlinear term can be estimated 
exactly in the same way as the original term $A^{1/2}(|u|^2)$ , and also the 
additional linear term is harmless. This remark was already used by \cite{GTV}.

Thus we have
\begin{theorem}
Let $k,l,b,b_1,p$ fulfill the assumptions of Theorem \ref{Theorem}. Let $u_0 
\in \widehat{H^{k,p}}$ , $n_0 \in \widehat{H^{l,p}}$ , $ n_1 \in 
\widehat{H^{l-1,p}}$ . Then the Cauchy problem (\ref{1}),(\ref{2}),(\ref{3}) is 
locally well-posed, i.e. there exists a unique solution
$$ u \in X^{k,b_1}_p(0,T) \, , \, n \in X^{l,b}_{+,p}(0,T) + X^{l,b}_{-,p}(0,T) 
\, , \, n_t \in X^{l-1,b}_{+,p}(0,T) + X^{l-1,b}_{-,p}(0,T) \, . $$
This solution satisfies
$$ u \in C^0([0,T],\widehat{H^{k,p}}) \, , \, n \in 
C^0([0,T],\widehat{H^{l,p}}) \, , \, n_t \in C^0([0,T],\widehat{H^{l-1,p}}) \, 
, $$
and the mapping data upon solution is locally Lipschitz continuous.
\end{theorem}
We use the notation $\langle \lambda \rangle := (1+\lambda^2)^{1/2}$ , and $a 
\pm$ to denote a number slightly larger (resp., smaller) than $a$ .
\section{Comparison with earlier results}
It is interesting to compare our results with those of \cite{GTV} for the case 
$p=2$. The lowest admissible choice in this case was $k=0$ , $l=-1/2$ , $p=2$. 
This is contained in our results, too.
\begin{itemize}
\item A choice, which improves this result from the scaling point of view for 
the Schr\"odinger part is $k=0$ , $p=1+\epsilon$ , $-\frac{2}{p'} < l \le - 
\frac{1}{p'} $ (with $b=b_1=\frac{1}{p}+$) and $\epsilon > 0$ small. It is 
easily checked that this choice is admissible due to Theorem \ref{Theorem}. \\
$\widehat{H^{k,p}}$ scales like $H^{\sigma}$ , where $\sigma = 
k-\frac{1}{p}+\frac{1}{2}$ , here: $\sigma = \frac{1}{2} - \frac{1}{1+\epsilon} 
\to -\frac{1}{2}$ ($\epsilon \to 0$),  $\widehat{H^{l,p}}$ scales like 
$H^{\lambda}$ , where $\lambda = l-\frac{1}{p}+\frac{1}{2}$ , here: $\lambda 
\to -\frac{1}{2}$ ($\epsilon \to 0$). \\
That $\widehat{H^{k,p}}$ scales like $H^{\sigma}$ here just means that (cf. 
(\ref{aa}) and (\ref{a})):
$$\| |\xi|^k \widehat{u_{\mu}}(\xi,0) 
\|_{L_{\xi}^{p'}} =:\|u_{\mu}(x,0)\|_{\widehat{\dot{H}^{k,p}}} = 
\mu^{k-\frac{1}{p}+\frac{3}{2}} 
\|u_0\|_{\widehat{\dot{H}^{k,p}}} $$
and 
$$ \|u_{\mu}(x,0)\|_{\dot{H}^{\sigma}} = \mu ^{\sigma +1} 
\|u_0\|_{\dot{H}^{\sigma}} $$
and the exponents of $\mu$ here coincide.
\item
Another admissible choice improving the result from the scaling point of view 
for the wave part is $k=0$ , $l=-\frac{1}{p} $ (with $b=b_1=\frac{1}{p}+$) and 
$2 \ge p > \frac{3}{2}.$ The conditions of Theorem \ref{Theorem} are fulfilled:
\begin{enumerate}
\item $k-l<2(1-b_1) \, \Leftrightarrow \, \frac{1}{p} <2(1-\frac{1}{p}) \, 
\Leftrightarrow \, p > \frac{3}{2}$ ,
\item $l \le 2k-\frac{1}{p'} \, \Leftrightarrow \, p \le 2 $ ,
\item $l+1-k < \frac{1}{p}+2(1-b) \, \Leftrightarrow \, -\frac{1}{p} +1 < 
\frac{1}{p} +2(1-\frac{1}{p}) $ , which is fulfilled, and
\item $l+1-k \le 2b_1 \, \Leftrightarrow \, -\frac{1}{p}+1 \le \frac{2}{p} \, 
\Leftrightarrow \, 1 \le \frac{3}{p} $ .
\end{enumerate}
$\widehat{H^{k,p}}$ scales like $H^{\sigma}$ with $\sigma = 
-\frac{1}{p}+\frac{1}{2} \to -\frac{1}{6}$ ($p \to \frac{3}{2}$) ,
$\widehat{H^{l,p}}$ scales like $H^{\lambda}$ with $\lambda = 
-\frac{2}{p}+\frac{1}{2} \to -\frac{5}{6}$ ($p \to \frac{3}{2}$) .
\end{itemize}
It is also interesting to remark that it is possible to choose $k<0$ and 
$l<-\frac{1}{2}$ (with a suitable $1<p<2$) , and nevertheless achieve local 
well-posedness for the Zakharov system (see details below). In this situation 
Holmer \cite{H} proved in the $L^2$-case that the mapping data upon solution in 
not $C^2$, so that a contraction mapping method as in our case cannot be 
applied. Moreover, the cubic nonlinear Schr\"odinger equation (\ref{NLS}) is 
known to be ill-posed for suitable data $u_0 \in \widehat{H^{k,p}}$ for any 
$-\frac{1}{p'} < k < 0 $ and $ p>1 $ (cf. \cite{G2}). This equation, as already 
remarked in the introduction, is the formal limit as $c \to \infty$ of a 
sequence of velocity-dependent Zakharov systems (replacing $\partial_x^2 - 
\partial_t^2$ by $c^{-2}\partial_x^2 - \partial_t^2$). So this limit must be 
singular in some sense.

In order to determine the minimal $k$, which fulfills all the assumptions in 
Theorem \ref{Theorem} we argue as follows: 
\begin{enumerate}
\item The conditions $2k>\frac{1}{p}-b_1$ and $k<l+2(1-b_1)$ require 
$\frac{1}{2p}-\frac{1}{2}b_1 < l+2-2b_1 \, \Leftrightarrow \, b_1 < 
\frac{2}{3}(l+2)-\frac{1}{3p}$ .
\item The conditions $2k \ge l+\frac{1}{p'}$ and $k<l+2(1-b_1)$ require 
$\frac{l}{2}+\frac{1}{2p'} < l+2-2b_1 \, \Leftrightarrow \, b_1 < 
\frac{l}{4}+\frac{3}{4}+\frac{1}{4p}$ .
\end{enumerate}
Thus $b_1$ has to be chosen such that $\frac{1}{p} < b_1 < 
\min(\frac{2}{3}(l+2)-\frac{1}{3p},\frac{l}{4}+\frac{3}{4}+\frac{1}{4p})$ , so 
that the condition $2k>\frac{1}{p}-b_1$ can only be fulfilled, if
\begin{equation}
\label{4.1}
2k  >  \frac{1}{p} - \frac{2}{3}(l+2) + \frac{1}{3p} = \frac{4}{3p} - 
\frac{2}{3}(l+2) 
\end{equation}
and
\begin{equation}
\label{4.2}
2k  >  \frac{1}{p} - \frac{l}{4} - \frac{3}{4} - \frac{1}{4p} = \frac{3}{4p} 
- \frac{1}{4}(l+3) \, .
\end{equation}
Moreover we need
\begin{equation}
\label{4.3}
2k \ge l + \frac{1}{p'} = l+1-\frac{1}{p} \, .
\end{equation} 
The lower bound for $2k$ in (\ref{4.1}) and (\ref{4.3}) is minimized, if
\begin{equation}
\label{4.4}
\frac{4}{3p}-\frac{2}{3}(l+2) = l+1-\frac{1}{p} \, \Leftrightarrow \, 
\frac{1}{p} = \frac{5}{7}l +1 \, .
\end{equation}
One easily checks that under this assumption all 3 lower bounds for $2k$ 
coincide. Thus we end up with (from (\ref{4.1})):
$$ 2k > \frac{4}{3}(\frac{5}{7}l+1)-\frac{2}{3}(l+2) = \frac{2}{7}l \quad 
\Leftrightarrow \quad k > \frac{l}{7} \, . $$
The minimal and optimal choice for $l$ here is $l=-\frac{1}{p}$ (because $l\ge 
-\frac{1}{p}$), which means by (\ref{4.4}): $p = \frac{12}{7}$ , and thus $ k > 
- \frac{1}{12} $ (from $k>\frac{l}{7}$) , and by (\ref{4.4}): 
$$ \frac{5}{7}l = \frac{1}{p}-1 = - \frac{5}{12} \, \Leftrightarrow \, l = - 
\frac{7}{12} \, . $$
Moreover, we should choose $b_1 < \frac{2}{3}(l+2)-\frac{1}{3p} = \frac{3}{4} 
$. \\
It is now completely elementary to see that that the choice 
$k=-\frac{1}{12}+\epsilon$ , $l=-\frac{7}{12},$  $b=b_1=\frac{3}{4}-\epsilon$ , 
$p=\frac{12}{7}$ ($\epsilon >0$ small)  meets all the assumptions of Theorem 
\ref{Theorem}. 

In this situation we have: $\widehat{H^{k,p}}$ scales like $H^{\sigma}$ with 
$\sigma = -\frac{1}{6} + \epsilon$ , and $\widehat{H^{l,p}}$ scales like 
$H^{\lambda}$ with $\lambda = -\frac{2}{3}$ . 

This is an improvement from the scaling point of view for both the 
Schr\"odinger and the wave part, compared to the $L^2$-result of \cite{GTV}, 
where $\sigma=0$ and $\lambda=-\frac{1}{2}.$ 
\section{Nonlinear estimates}
In order to estimate the nonlinearities we use the following simple application 
of H\"older's inequality.
\begin{lemma}
\label{Lemma 1.1}
For $1/p + 1/p' = 1$ , $1<p<\infty$ , the following estimate holds:
\begin{eqnarray*}
\lefteqn{ | \int\int \widehat{v}(\zeta) \widehat{v_1}(\zeta_1) 
\widehat{v_2}(\zeta_2) K(\zeta_1,\zeta_2) d\zeta_1 d\zeta_2 | } \\
& \le  & \sup_{\zeta_1} \left( \int |K(\zeta_1,\zeta_2)|^p d\zeta_2 
\right)^{1/p} \|\widehat{v_1}\|_{L^p} \|\widehat{v}\|_{L^{p'}} 
\|\widehat{v_2}\|_{L^{p'}} \, ,
\end{eqnarray*}
where $\zeta := \zeta_1 - \zeta_2$ .
\end{lemma}
{\bf Proof:} 
\begin{eqnarray*}
\lefteqn{  | \int\int \widehat{v}(\zeta) \widehat{v_1}(\zeta_1) 
\widehat{v_2}(\zeta_2) K(\zeta_1,\zeta_2) d\zeta_1 d\zeta_2 | }\\
 & \le & \|\widehat{v_1}\|_{L^p} (\int|\int \widehat{v}(\zeta_1 - \zeta_2) 
\widehat{v_2}(\zeta_2) K(\zeta_1,\zeta_2) d\zeta_2|^{p'} d\zeta_1 
)^{1/p'} \\
& \le & \|\widehat{v_1}\|_{L^p} \{\int[(\int|\widehat{v}(\zeta_1 - \zeta_2) 
\widehat{v_2}(\zeta_2)|^{p'} d\zeta_2)(\int |K(\zeta_1,\zeta_2)|^p 
d\zeta_2)^{\frac{p'}{p}}] d\zeta_1 \}^{1/p'} \\ 
& \le & \|\widehat{v_1}\|_{L^p} (\sup_{\zeta_1} \int|K(\zeta_1,\zeta_2)|^p 
d\zeta_2)^{1/p} (\int\int |\widehat{v}(\zeta_1 - \zeta_2) 
\widehat{v_2}(\zeta_2)|^{p'} d\zeta_2 d\zeta_1 )^{1/p'} \\
& \le & \sup_{\zeta_1} \left(\int|K(\zeta_1,\zeta_2)|^p d\zeta_2\right)^{1/p} 
\|\widehat{v_1}\|_{L^p} \|\widehat{v}\|_{L^{p'}} \|\widehat{v_2}\|_{L^{p'}} \, 
.
\end{eqnarray*}
{\bf Remark:} Similarly one can prove
\begin{eqnarray*}
\lefteqn{ | \int\int \widehat{v}(\zeta) \widehat{v_1}(\zeta_1) 
\widehat{v_2}(\zeta_2) K(\zeta_1,\zeta_2) d\zeta_1 d\zeta_2 | } \\
& \le  & \sup_{\zeta} \left( \int |K(\zeta + \zeta_2,\zeta_2)|^p d\zeta_2 
\right)^{1/p} \|\widehat{v}\|_{L^p} \|\widehat{v_1}\|_{L^{p'}} 
\|\widehat{v_2}\|_{L^{p'}} \, .
\end{eqnarray*}

Our first aim is to estimate the nonlinearity $ f = n_{\pm} u $ in 
$X^{k,-c_1}_p$ for given $n_{\pm} \in X^{l,b}_{\pm,p}$ and $u \in X^{k,b_1}_p$ 
. We estimate $\widehat{f}(\xi_1',\tau_1) = (\widehat{n_{\pm}} * 
\widehat{u})(\xi_1',\tau_1)$ in terms of $\widehat{n_{\pm}}(\xi,\tau)$ and 
$\widehat{u}(\xi_2',\tau_2)$ , where $\xi = \xi_1' - \xi_2'$ , $ \tau = \tau_1 
- \tau_2$ . We also introduce the variables $\sigma_1 = \tau_1 + \xi_1'^2$ , 
$\sigma_2 = \tau_2 + \xi_2'^2$ , $ \sigma = \tau \pm |\xi| $ , so that
\begin{equation}
\label{**}
z := \xi_1'^2 - \xi_2'^2 \mp |\xi| = \sigma_1 - \sigma_2 - \sigma \, .
\end{equation}
Define $\widehat{v_2} = \langle \xi_2' \rangle^k \langle \sigma_2 \rangle^{b_1} 
\widehat{u} $ and $ \widehat{v} = \langle \xi \rangle^l \langle \sigma 
\rangle^b \widehat{n_{\pm}} $ , so that $ \|u\|_{X^{k,b_1}_p} = \| 
\widehat{v_2}\|_{L^{p'}}$ and $ \|n_{\pm}\|_{X^{l,b}_{\pm,p}} = 
\|\widehat{v}\|_{L^{p'}}$ . In order to estimate $f$ in $X^{k,-c_1}_p$ we take 
its scalar product with a function in $X^{-k,c_1}_{p'}$ with Fourier transform 
$\langle \xi_1' \rangle^k \langle \sigma_1 \rangle^{-c_1} \widehat{v_1}$ with 
$\widehat{v_1} \in L^p$ .

In the sequel we want to show an estimate of the form
$$ |S| \le c \|\widehat{v}\|_{L^{p'}} \|\widehat{v_1}\|_{L^p} 
\|\widehat{v_2}\|_{L^{p'}} \, , $$
where
$$ S:= \int \frac{|\widehat{v} \widehat{v_1} \widehat{v_2}| \langle 
\xi_1'\rangle^k}{\langle \sigma \rangle^b \langle \sigma_1 \rangle^{c_1} 
\langle \sigma_2 \rangle^{b_1} \langle \xi_2' \rangle^k \langle \xi \rangle^l} 
d\xi_1' d\xi_2' d\tau_1 d\tau_2 \, . $$
This directly gives the desired estimate
\begin{equation}
\label{*}
\|n_{\pm} u \|_{X^{k,-c_1}_p} \le c \|n_{\pm}\|_{X^{l,b}_{\pm,p}} 
\|u\|_{X^{k,b_1}_p} \, .
\end{equation}
\begin{prop}
\label{Proposition 1.1}
The estimate (\ref{*}) holds under the following assumptions:
$$ k \ge 0 \, , \, l \ge -1/p \, , \, k-l \le 2c_1 \, , \, k-l \le 2/p \, , $$
where $ c_1 \ge 0 \, , \, b > 1/p \, , \, b_1 > 1/p \, , \, 1 < p \le 2 $ .
\end{prop}
{\bf Remark:} We simplify (\ref{**}) as follows. If (\ref{**}) holds with the 
minus sign and if $\xi_1' \ge \xi_2'$ (resp., $\xi_1' \le \xi_2'$) , we have
$$ z = \xi_1'^2 - \xi_2'^2 - |\xi_1' - \xi_2'| = (\xi_1' \mp 1/2)^2 - (\xi_2' 
\mp 1/2)^2 = \xi_1^2 - \xi_2^2 \, $$
where $ \xi_i = \xi_i' \mp 1/2 $ . Thus the region $\xi_1' \ge \xi_2'$ (resp., 
$\xi_1' \le \xi_2'$) of $S$ is majorized by
$$ \overline{S} = c \int \frac{|\widehat{v}(\xi,\tau) \widehat{v_1}(\xi_1 \pm 
1/2,\tau_1) \widehat{v_2}(\xi_2 \pm 1/2,\tau_2)| \langle 
\xi_1\rangle^k}{\langle \sigma \rangle^b \langle \sigma_1 \rangle^{c_1} \langle 
\sigma_2 \rangle^{b_1} \langle \xi_2 \rangle^k \langle \xi \rangle^l} d\xi_1 
d\xi_2 d\tau_1 d\tau_2 \, , $$
where now
\begin{eqnarray}
\label{***}
z = \xi_1^2 - \xi_2^2 = \sigma_1 - \sigma_2 - \sigma \quad , \quad \xi = \xi_1 
- \xi_2 \quad , \quad \tau = \tau_1 - \tau_2 \\ \nonumber
\sigma_i = \tau_i + (\xi_i \pm 1/2)^2 \quad , \quad \sigma = \tau \pm |\xi| = 
\tau \pm |\xi_1 - \xi_2| \, .
\end{eqnarray}
Also, the plus sign in (\ref{**}) can be treated similarly by again defining 
$\xi_i = \xi_i' \pm 1/2$ . If one wants to estimate $\overline{S}$ by $c 
\|\widehat{v}\|_{L^{p'}} \|\widehat{v_1}\|_{L^p} \|\widehat{v_2}\|_{L^{p'}}$ , 
the variables $\xi_i$ and $\xi_i \pm 1/2$ are completely equivalent, thus we do 
not distinguish between them.\\
{\bf Proof of Proposition \ref{Proposition 1.1}:}
According to Lemma \ref{Lemma 1.1} we have to show
$$ C^p := \sup_{\xi_1,\sigma_1} \langle \sigma_1 \rangle^{-c_1p} \langle \xi_1 
\rangle^{kp} \int \frac{d\xi_2 d\sigma_2}{\langle \sigma \rangle^{bp} \langle 
\sigma_2 \rangle^{b_1p} \langle \xi \rangle^{lp} \langle \xi_2 \rangle^{kp}} < 
\infty \, . $$
{\bf Case 1:} $|\xi_1| \le 2 |\xi_2| $ ($\Rightarrow \, |\xi| \le 3|\xi_2|$) . 
\\
If $|\xi_2| \le 1$ we have $ \langle \xi \rangle \sim 1 $ and thus
$$ C^p \le c \sup_{\xi_1,\sigma_1} \int_{|\xi_2| \le 1} d\xi_2 \int_0^{\infty} 
\langle \sigma_2 \rangle^{-b_1 p} d\sigma_2 < \infty \, , $$
because $b_1 > 1/p$ . If $|\xi_2| \ge 1$ we get from (\ref{***}) for 
$\xi_1,\sigma_1,\sigma_2$ fixed: $ \frac{d\sigma}{d\xi_2} = 2\xi_2 $ , and thus
$$ C^p \le c \sup_{\xi_1,\sigma_1} \int \langle \xi \rangle^{-lp} \langle \xi_2 
\rangle^{-1} \langle \sigma \rangle^{-bp} \langle \sigma_2 \rangle^{-b_1p} 
d\sigma d\sigma_2 \, . $$
For $ l \ge 0$ we immediately have
$$ C^p \le c \int \langle \sigma \rangle^{-bp} d\sigma \int \langle \sigma_2 
\rangle^{-b_1 p} d\sigma_2 < \infty \, , $$
whereas for $l \le 0$ we use our assumption $l \ge -1/p$ and get the same bound 
by using $ \langle \xi \rangle^{-lp} \langle \xi_2 \rangle^{-1} \le c \langle 
\xi_2 \rangle^{-lp-1} \le c $ .\\
{\bf Case 2:} $|\xi_1| \ge 2 |\xi_2| $ ($\Rightarrow \, |\xi| \sim |\xi_1|$) . 
\\
From (\ref{***}) we conclude $\xi_1^2 \le c(|\sigma_1| + |\sigma_2| + |\sigma|) 
$ and distinguish three cases. \\
{\bf Case 2a:} $|\sigma_1|$ dominant, i.e. $|\sigma_1| \ge |\sigma_2|,|\sigma|$ 
, ($\Rightarrow \, \xi_1^2 \le c|\sigma_1|$) . \\
This implies, using our assumption $k-l \le 2c_1$ :
\begin{eqnarray*}
C^p & \le & c \sup_{\xi_1,\sigma_1} \langle \xi_1 \rangle^{(k-l-2c_1)p} \int 
d\xi_2 
d\sigma_2 \langle \sigma \rangle^{-bp} \langle \sigma_2 \rangle^{-b_1p} \\
&\le & c \sup_{\xi_1,\sigma_1}\int d\xi_2 d\sigma_2 \langle \sigma 
\rangle^{-bp} \langle \sigma_2 
\rangle^{-b_1 p} \, .
\end{eqnarray*}
If $|\xi_2| \le 1$ this is immediately bounded by $c\int_{|\xi_2|\le 1} d\xi_2 
\int d\sigma_2 \langle \sigma_2 \rangle^{-b_1 p} < \infty $ , whereas for 
$|\xi_2| \ge 1$ we use $\frac{d\sigma}{d\xi_2} = 2\xi_2 \sim 2\langle \xi_2 
\rangle$ again and get the bound 
$$c \int d\sigma \langle \sigma \rangle^{-bp} \int d\sigma_2 \langle \sigma_2 
\rangle^{-b_1 p} < \infty \, ,$$ 
using $b,b_1 > 1/p$ . \\
{\bf Case 2b:} $|\sigma_2|$ dominant ( $\Rightarrow \, \xi_1^2 \le c 
|\sigma_2|$ ).\\
We have
$$ C^p \le c \sup_{\xi_1,\sigma_1} \langle \xi_1 \rangle^{(k-l)p} \int d\xi_2 
d\sigma_2 \langle \sigma \rangle^{-bp} \langle \sigma_2 \rangle^{-b_1p} \langle 
\xi_2 \rangle^{-kp} \, . $$
The case $k \le l$ is simple and leads to
$$ C^p \le c \sup_{\xi_1,\sigma_1} \int d\xi_2 d\sigma_2 \langle \sigma 
\rangle^{-bp} \langle \sigma_2 \rangle^{-b_1p} \, , $$
which can be handled as in case 2a.\\
The case $k>l$ is treated as follows:
$$ C^p \le c \sup_{\xi_1,\sigma_1} \int d\xi_2 d\sigma_2 \langle \sigma 
\rangle^{-bp} \langle \sigma_2 \rangle^{-b_1p+\frac{(k-l)p}{2}} \langle \xi_2 
\rangle^{-kp} \, . $$
Substituting $ y = \xi_2^2$ , thus $d\xi_2 = \frac{dy}{2|y|^{1/2}}$ , leads to
$$ C^p \le c \sup_{\xi_1,\sigma_1} \int |y|^{-\frac{1}{2}} \langle y 
\rangle^{-\frac{kp}{2}} \int \langle \sigma \rangle^{-bp} \langle \sigma_2 
\rangle^{-b_1p+\frac{(k-l)p}{2}} d\sigma_2  dy \, . $$
From (\ref{***}) we have $\langle \sigma \rangle = \langle \sigma_2 
-(\sigma_1-\xi_1^2 +y) \rangle $ , and thus by \cite{GTV}, Lemma 4.2, using 
$b,b_1 > 1/p$ :
$$ \int d\sigma_2 \langle \sigma \rangle^{-bp} \langle \sigma_2 \rangle^{-b_1p 
+ \frac{(k-l)p}{2}} \le c \langle \sigma_1 - \xi_1^2 +y 
\rangle^{-1+\frac{(k-l)p}{2} - } \, , $$
because $-b_1p + \frac{(k-l)p}{2} < -1+1 = 0$ using our assumptions $b_1 > 1/p$ 
and $k-l \le 2/p \, .$ Thus 
$$ C^p \le c \sup_{\xi_1,\sigma_1} \int |y|^{-\frac{1}{2}} \langle y 
\rangle^{-\frac{kp}{2}} \langle \sigma_1 - \xi_1^2 +y 
\rangle^{-1+\frac{(k-l)p}{2}-} dy \, . $$
The supremum occurs for $\sigma_1 = \xi_1^2$ by \cite{GTV}, Lemma 4.1, so that
$$ C^p  \le  c \int |y|^{-\frac{1}{2}} \langle y \rangle^{-\frac{kp}{2}} 
\langle y \rangle^{-1+\frac{(k-l)p}{2}-} dy 
  =  \int |y|^{-\frac{1}{2}} \langle y \rangle^{-1-\frac{lp}{2}-} dy < \infty 
\, , $$
because $l \ge -1/p$ . \\
{\bf Case 2c:} $|\sigma|$ dominant.\\
This case can be treated like case 2b, which completes the proof.\\[1cm]
It is also possible to prove (\ref{*}) in certain cases where $k$ is negative. 
This is done in the following
\begin{prop}
\label{Proposition 1.2}
Estimate (\ref{*}) holds under the following conditions:
$$ k \le 0 \, , \, l \ge -\frac{1}{p} \, , \,  k \ge -\frac{1}{p} \, , $$
$$k-l \le 2c_1 \, , \, k-l \le 
\frac{2}{p} \, , \, l+k > \frac{1}{p} -2b_1 \, , \, l+k > \frac{1}{p} -2b \, , 
\, l+k \ge -\frac{1}{p} -2c_1  \, , $$
where $ c_1 \ge 0 $ and $b_1,b > \frac{1}{p}$ .
\end{prop}
{\bf Proof:} \\
{\bf Case 1:} $|\xi_1| \sim |\xi_2|$ \\
This case can be treated exactly like case 1 in the previous proposition, using 
$l \ge -1/p$ . \\
{\bf Case 2:} $|\xi_1| << |\xi_2| $ ( $\Rightarrow \, |\xi| \sim |\xi_2|$ ).\\
Using the notation of the previous proposition we have
$$ C^p \le c \sup_{\xi_1,\sigma_1} \langle \sigma_1 \rangle^{-c_1p} \langle 
\xi_1 \rangle^{kp} \int d\xi_2 d\sigma_2 \langle \xi_2 \rangle^{-lp-kp} \langle 
\sigma \rangle^{-bp} \langle \sigma_2 \rangle^{-b_1 p} \, . $$
From (\ref{***}) we get $\xi_2^2 \le c(|\sigma_1|+|\sigma_2|+|\sigma|)$ . \\
{\bf Case 2a:} $|\sigma_1|$ dominant ( $ \Rightarrow \, \xi_2^2 \le c 
|\sigma_1|$ ) . \\
Thus
$$ C^p \le c \sup_{\xi_1,\sigma_1} \int d\xi_2 d\sigma_2 \langle \xi_2 
\rangle^{-lp-kp-2c_1p} \langle \sigma \rangle^{-bp} \langle \sigma_2 
\rangle^{-b_1 p} \, . $$ 
If $|\xi_2| \le 1$ we easily get the bound
$$ C^p \le c \int_{|\xi_2| \le 1} d\xi_2 \int \langle \sigma_2 \rangle^{-b_1p} 
d\sigma_2 < \infty \, . $$
If $|\xi_2| \ge 1$ we use (\ref{***}) and get for fixed 
$\xi_1,\sigma_1,\sigma_2:$ $\frac{d\sigma}{d\xi_2} = 2|\xi_2| \sim 2\langle 
\xi_2 \rangle$ , so that, using the condition $l+k \ge -\frac{1}{p} -2c_1$ , we 
have $ \langle \xi_2 \rangle^{-(l+k+2c_1)p-1} \le c $ , and thus the bound
$$ C^p \le c \int d\sigma \langle \sigma \rangle^{-bp} \int d\sigma_2 \langle 
\sigma_2 \rangle^{-b_1 p} < \infty $$
by $b,b_1 > 1/p$ .\\
{\bf Case 2b:} $|\sigma_2|$ dominant ($\Rightarrow \xi_2^2 \le c|\sigma_2|$ ) . 
\\
Ignoring the factor $\langle \sigma_1 \rangle^{-c_1 p} \langle \xi_1 
\rangle^{kp}$ we estimate
$$ C^p \le c \sup_{\xi_1,\sigma_1} \int d\xi_2 d\sigma_2 \langle \xi_2 
\rangle^{-lp-kp+2-2b_1 p+} \langle \sigma \rangle^{-bp} \langle \sigma_2 
\rangle^{-1-} \, . $$
The case $|\xi_2| \le 1$ is easy again, and for $|\xi_2| \ge 1$ we again use 
$\frac{d\sigma}{d\xi_2} =2|\xi_2| \sim 2\langle \xi_2 \rangle$ and $\langle 
\xi_2 \rangle^{-lp-kp+2-2b_1 p-1+} \le c$ , using our assumption 
$l+k>\frac{1}{p} -2b_1$ , and thus
$$ C^p \le c \sup_{\xi_1,\sigma_1} \int d\sigma \langle \sigma \rangle^{-bp} 
\int d\sigma_2 \langle \sigma_2 \rangle^{-1-} < \infty \, . $$ \\
{\bf Case 2c:} $|\sigma|$ dominant . \\
This case can be handled like case 2b, using the assumption 
$l+k>\frac{1}{p}-2b.$ \\
{\bf Case 3:} $|\xi_1| >> |\xi_2|$ ( $\Rightarrow \, |\xi| \sim |\xi_1|$ ). \\
We have
$$ C^p \le c \sup_{\xi_1,\sigma_1} \langle \sigma_1 \rangle^{-c_1 p} \langle 
\xi_1 \rangle^{(k-l)p} \int d\xi_2 d\sigma_2 \langle \xi_2 \rangle^{-kp} 
\langle \sigma \rangle^{-bp} \langle \sigma_2 \rangle^{-b_1 p} \, . $$
Again we distinguish three cases. \\
{\bf Case 3a:} $|\sigma_1|$ dominant ( $\Rightarrow \, \xi_1^2 \le c|\sigma_1|$ 
) . \\
By our assumption $k-l \le 2c_1$ we get
\begin{eqnarray*}
C^p & \le & c \sup_{\xi_1,\sigma_1} \langle \xi_1 \rangle^{(k-l-2c_1)p} \int 
d\xi_2 d\sigma_2 \langle \xi_2 \rangle^{-kp} \langle \sigma \rangle^{-bp} 
\langle \sigma_2 \rangle^{-b_1 p} \\
& \le & c \sup_{\xi_1,\sigma_1} \int d\xi_2 d\sigma_2 \langle \xi_2 
\rangle^{-kp} \langle \sigma \rangle^{-bp} \langle \sigma_2 \rangle^{-b_1 p} \, 
,
\end{eqnarray*}
which is easily handled for $|\xi_2| \le 1$ , whereas for $|\xi_2| \ge 1$ , 
using again $\frac{d\sigma}{d\xi_2} = 2|\xi_2| \sim 2\langle \xi_2 \rangle $ 
and our assumption $ k \ge -1/p$ , we arrive at
$$ C^p \le c \int\langle \sigma \rangle^{-bp} d\sigma \int \langle \sigma_2 
\rangle^{-b_1 p} d\sigma_2 < \infty \, . $$ 
{\bf Case 3b:} $|\sigma_2|$ dominant ( $\Rightarrow \, \xi_1^2 \le c 
|\sigma_2|$ ). \\
$$ C^p \le c \sup_{\xi_1,\sigma_1} \langle \xi_1 \rangle^{(k-l)p} \int d\xi_2 
d\sigma_2 \langle \xi_2 \rangle^{-kp} \langle \sigma \rangle^{-bp} \langle 
\sigma_2 \rangle^{-b_1 p} \, . $$
In the case $k\le l$ we have
$$ C^p \le c \sup_{\xi_1,\sigma_1} \int d\xi_2 d\sigma_2 \langle \xi_2 
\rangle^{-kp} \langle \sigma \rangle^{-bp} \langle \sigma_2 \rangle^{-b_1 p} \, 
, $$
which is simple to handle for $|\xi_2| \le 1$ , and for $|\xi_2| \ge 1$ we use 
$\frac{d \sigma}{d\xi_2} = 2|\xi_2| \sim 2 \langle \xi_2 \rangle$ and the 
assumption $k \ge -1/p$ and estimate
$$ C^p \le c \int d\sigma \langle \sigma \rangle^{-bp} \int d\sigma_2 \langle 
\sigma_2 \rangle^{-b_1 p} < \infty \, . $$
In the case $k \ge l$ we use $\xi_1^2 \le c |\sigma_2|$ and the substitution 
$y=\xi_2^2$ and $d\xi_2 = \frac{dy}{2|y|^{1/2}}$ to get
\begin{eqnarray*}
C^p & \le & c \sup_{\xi_1,\sigma_1} \int d\xi_2 d\sigma_2 \langle \xi_2 
\rangle^{-kp} \langle \sigma \rangle^{-bp} \langle \sigma_2 \rangle^{-b_1 
p+\frac{(k-l)p}{2}} \\
& \le & c \sup_{\xi_1,\sigma_1} \int |y|^{-\frac{1}{2}} \langle y 
\rangle^{-\frac{kp}{2}} \int  \langle \sigma \rangle^{-bp} \langle \sigma_2 
\rangle^{-b_1 p+\frac{(k-l)p}{2}} d\sigma_2 dy \, .
\end{eqnarray*}
Now $\langle \sigma \rangle = \langle \sigma_2 - (\sigma_1 - \xi_1^2 +y) 
\rangle $ by (\ref{***}) , and \cite{GTV}, Lemma 4.2 implies
$$ \int \langle \sigma \rangle^{-bp} \langle \sigma_2 \rangle^{b_1 
p+\frac{(k-l)p}{2}} d\sigma_2 \le c \langle \sigma_1 - \xi_1^2 +y 
\rangle^{-1+\frac{(k-l)p}{2}-} \, ,$$
where we used the assumption $k-l \le 2/p$ . Thus by \cite{GTV} , Lemma 4.3 , 
using $k \ge -1/p$ and $l \ge -1/p$ :
\begin{eqnarray*}
C^p & \le & c\sup_{\xi_1,\sigma_1} \int |y|^{-\frac{1}{2}} \langle y 
\rangle^{-\frac{kp}{2}} \langle \sigma_1 - \xi_1^2 +y 
\rangle^{-1+\frac{(k-l)p}{2}-} dy \\
& \le & c\int |y|^{-\frac{1}{2}} \langle y \rangle^{-\frac{kp}{2}} \langle y 
\rangle^{-1+\frac{(k-l)p}{2}-} dy \\
& = & c\int |y|^{-\frac{1}{2}} \langle y \rangle^{-1-\frac{lp}{2}-} dy \, < 
\infty \, .
\end{eqnarray*}
{\bf Case 3c:} $|\sigma|$ dominant. \\
This case is treated like Case 3b. \\[1cm]
Next, we want to estimate the nonlinearity $(|u|^2)_x$, namely
$$ |W| \le c \|\widehat{v}\|_{L^p} \|\widehat{v_1}\|_{L^{p'}} 
\|\widehat{v_2}\|_{L^{p'}} \, , $$
where
$$ W:= \int \frac{|\widehat{v} \widehat{v_1} \widehat{v_2}| \langle \xi 
\rangle^l |\xi|}{\langle \sigma \rangle^c \langle \sigma_1 \rangle^{b_1} 
\langle \sigma_2 \rangle^{b_1} \langle \xi_1' \rangle^k \langle \xi_2' 
\rangle^k} d\xi_1' d\xi_2' d\tau_1 d\tau_2 \, . $$
This implies the desired estimate
\begin{equation}
\label{W}
\| (|u|^2)_x\|_{X^{l,-c}_{\pm,p}} \le c \|u\|^2_{X^{k,b_1}_p} \, . 
\end{equation}
\begin{prop}
\label{Proposition 1.3}
Estimate (\ref{W}) holds under the following conditions:
$$ k \ge 0 \, , \, l \le 2k-\frac{1}{p'} \, , \, l+1-k \le \frac{1}{p} + 2c \, 
, \, l+1-k \le 2b_1 \, , $$
where $ c \ge 0 \, , \, b_1 > 1/p \, , \, 1 < p \le 2 $ .
\end{prop}
{\bf Proof:} According to the remark after Lemma \ref{Lemma 1.1} we have to 
show
\begin{equation}
\label{4}
C^p := \sup_{\xi,\sigma} \langle \sigma \rangle^{-cp} \langle \xi \rangle^{lp} 
|\xi|^p \int d\xi_2 d\sigma_2 \langle \xi_1 \rangle^{-kp} \langle \xi_2 
\rangle^{-kp} \langle \sigma_1 \rangle^{-b_1 p} \langle \sigma _2 \rangle^{-b_1 
p} < \infty \, .
\end{equation}
{\bf Case 1:} $|\xi_1| \sim |\xi_2|$ ( $\Rightarrow |\xi| \le 
c|\xi_1|,c|\xi_2|$ ) . \\
Applying the remark after Proposition \ref{Proposition 1.1} we have
$$ \xi_1^2 - \xi_2^2 = (\xi + \xi_2)^2 - \xi_2^2 = \sigma_1 - \sigma_2 -  
\sigma \, . $$
Thus, for fixed $\xi,\sigma,\sigma_2$, we have $\frac{d\sigma_1}{d\xi_2} = 
2\xi$ , so that
$$ C^p \le c \sup_{\xi,\sigma} \langle \xi \rangle^{lp} |\xi|^{p-1} \langle \xi 
\rangle^{-2kp} \int \langle \sigma_1 \rangle^{-b_1 p} d\sigma_1 \int \langle 
\sigma_2 \rangle^{-b_1 p} d\sigma_2 \, . $$
This is easily seen to be finite under the assumption $l-2k+1 \le 1/p $ $ 
\Leftrightarrow$ $ l \le 2k - \frac{1}{p'}. $  \\
{\bf Case 2:} $|\xi_1| >> |\xi_2|$ ( $\Rightarrow \, |\xi| \sim |\xi_1|$ ) (and 
analogously $|\xi_2| >> |\xi_1|$ ) . \\
{\bf Case 2a:} $|\sigma|$ dominant ( $\Rightarrow \, \xi^2 \sim \xi_1^2 \le 
c|\sigma|$ ) . \\
Using the relation $\frac{d\sigma_1}{d\xi_2} = 2\xi$ again and ignoring the 
term $\langle \xi_2 \rangle^{-kp}$ we arrive at
\begin{eqnarray*}
C^p & \le & c \sup_{\xi,\sigma} \langle \xi \rangle^{(-2c+l-k)p} |\xi|^p \int 
d\xi_2 d\sigma_2 \langle \sigma_1 \rangle^{-b_1 p} \langle \sigma _2 
\rangle^{-b_1 p} \\
& \le & c \sup_{\xi,\sigma} \langle \xi \rangle^{(-2c+l-k)p} |\xi|^{p-1} \int 
\langle \sigma_1 \rangle^{-b_1 p} d\sigma_1 \int \langle \sigma _2 
\rangle^{-b_1 p} d\sigma_2 \, ,
\end{eqnarray*}
which can be seen to be finite under the assumption $l+1-k \le \frac{1}{p} +2c$ 
. \\
{\bf Case 2b:} $|\sigma_1|$ dominant ( $\Rightarrow \, \xi^2 \sim \xi_1^2 \le c 
|\sigma_1| $ ) (and similarly $|\sigma_2|$ dominant).\\
We have
$$ C^p \le c \sup_{\xi,\sigma} \langle \xi \rangle^{(l-k)p} |\xi|^p \int d\xi_2 
d\sigma_2 \langle \xi_2 \rangle^{-kp} \langle \sigma_1 \rangle^{-b_1 p} \langle 
\sigma _2 \rangle^{-b_1 p} \, . $$
{\bf Case 2b$\alpha$:} In the case $l+1-k \ge 1/p$ we introduce the variable 
$z:= \xi_1^2 - \xi_2^2 = (\xi +\xi_2)^2 - \xi_2^2$ and get for fixed $\xi$: 
$\frac{dz}{d\xi_2} = 2 \xi$ . We also have $z=\xi^2 +2\xi \xi_2 $ $ 
\Leftrightarrow$ $ \xi_2 = \frac{z-\xi^2}{2\xi}$ and $z \sim \xi_1^2 \sim 
\xi^2$ , so that we get
\begin{eqnarray*}
C^p & \le & c \sup_{\xi,\sigma} \langle \xi \rangle^{(l-k)p} \langle \xi 
\rangle^{p-1} \int_{0\le z \le c\xi^2} dz d\sigma_2 \langle \sigma_1 
\rangle^{-b_1 p} \langle \sigma_2 \rangle^{-b_1 p} \langle \frac{z-\xi^2}{2\xi} 
\rangle^{-kp} \\
& \le & c \sup_{\xi,\sigma} \langle \xi \rangle^{-1} \int_{0\le z \le c\xi^2} 
dz d\sigma_2 \langle \sigma_1 \rangle^{\frac{(l-k+1)p}{2}-b_1 p} \langle 
\sigma_2 \rangle^{-b_1 p} \langle \frac{z-\xi^2}{2\xi} \rangle^{-kp} \, .
\end{eqnarray*}
Now we have by \cite{GTV}, Lemma 4.2:
$$ \int d\sigma_2 \langle \sigma_1 \rangle^{\frac{(l-k+1)p}{2}-b_1 p} \langle 
\sigma_2 \rangle^{-b_1 p} \le c \langle \sigma + 
z\rangle^{\frac{(l-k+1)p}{2}-b_1 p} \, , $$
where we used the assumption $l-k+1 \le 2b_1$ as well as (\ref{***}), namely 
$\sigma_1 - \sigma_2 = \sigma + z,$ so that we arrive at
$$ C^p \le c \sup_{\xi,\sigma} \langle \xi \rangle^{-1} \int_{0}^{c\xi^2} dz 
\langle \sigma + z \rangle^{\frac{(l-k+1)p}{2}-b_1 p} \langle 
\frac{z-\xi^2}{2\xi} \rangle^{-kp} \, . $$
With $y:= z-\xi^2$ we have
\begin{eqnarray*}
C^p & \le &  c \sup_{\xi,\sigma} \langle \xi \rangle^{-1} 
\int_{-c\xi^2}^{c\xi^2} \langle \sigma + \xi^2 +y 
\rangle^{\frac{(l-k+1)p}{2}-b_1 p} \langle \frac{y}{2\xi} \rangle^{-kp} dy \\
& = & c \sup_{\xi} \langle \xi \rangle^{-1} \int_{-c\xi^2}^{c\xi^2} \langle y 
\rangle^{\frac{(l-k+1)p}{2}-b_1 p} \langle \frac{y}{2\xi} \rangle^{-kp} dy \, ,
\end{eqnarray*}
where we used \cite{GTV}, Lemma 4.3. The case $|\xi| \le 1$ is easily handled. 
If $|\xi| \ge 1$ we get \\
{\bf a:} in the region $|y| \le |\xi|$ we have $\langle \frac{y}{2\xi} 
\rangle^{-kp} \sim 1$ and  $ \langle y \rangle^{\frac{(l-k+1)p}{2}-b_1p} \le 1$ 
by our assumption $l-k+1 \le 2b_1$ , so that the integral is bounded by 
$c|\xi|$ , thus $C^p < \infty.$ \\
{\bf b:} In the region $|\xi| \le |y| \le c\xi^2$ we have $\langle 
\frac{y}{2\xi}\rangle \sim |\frac{y}{2\xi}| \sim \frac{\langle y 
\rangle}{\langle \xi \rangle}$ , so that
$$ C^p \le c \sup_{\xi} \langle \xi \rangle^{-1} \langle \xi \rangle^{kp} 
\int_{|\xi|}^{c\xi^2} \langle y \rangle^{\frac{(l-k+1)p}{2}-b_1 p -kp} dy \, . 
$$ 
If $ \frac{(l-k+1)p}{2} -b_1 p -kp < -1 $ we have
$$ C^p \le c \sup_{\xi} \langle \xi 
\rangle^{-1+kp+\frac{(l-k+1)p}{2}-b_1p-kp+1} = c \sup_{\xi} \langle \xi 
\rangle^{\frac{(l-k+1)p}{2}-b_1p} \, , $$
which is finite under the assumption $l-k+1 \le 2b_1$ .\\
If $ \frac{(l-k+1)p}{2} -b_1 p -kp \ge -1 $ we have
$$ C^p \le c \sup_{\xi} \langle \xi \rangle^{-1+kp+(l-k+1)p-2b_1p-2kp+2+} = c 
\sup_{\xi} \langle \xi \rangle^{(l-2k+1)p-2b_1p+1+} < \infty \, , $$
because $(l-2k+1)p-2b_1p+1+ \le 0$ $\Leftrightarrow$ $ l-2k < 
2b_1-1-\frac{1}{p}$ , which is fulfilled under the assumption $l-2k \le 
\frac{1}{p} - 1 = - \frac{1}{p'}$ for $b_1 > \frac{1}{p}$ .\\
{\bf Case 2b$\beta$:}
In the case $l+1-k \le \frac{1}{p}$ we directly get by 
$\frac{d\sigma_1}{d\xi_2} = 2\xi$ (for $\sigma,\xi,\sigma_2$ fixed), ignoring 
the term $\langle \xi_2 \rangle^{-kp}$ in the integral:
$$ C^p \le c \sup_{\xi,\sigma} \langle \xi \rangle^{(l-k)p} |\xi|^{p-1} \int 
\langle \sigma_1 \rangle^{-b_1 p} d\sigma_1 \int \langle \sigma_2 \rangle^{-b_1 
p} d\sigma_2 < \infty \, . $$

The case $k<0$ is considered in the following
\begin{prop}
\label{Proposition 1.4}
Estimate (\ref{W}) holds under the following conditions:
$$ k \le 0 \, , \, l \le 2k-\frac{1}{p'} \, , \, 2k \ge -c \, , \, 2k > 
\frac{1}{p} - b_1 \, , $$
where $ \frac{1}{p'} > c \ge 0 \, , \, b_1 > \frac{1}{p} \, , \, 1 < p \le 2$ .
\end{prop}
{\bf Proof:} We again have to show (\ref{4}) as in the previous proof.\\
{\bf Case 1:} $|\xi_1| \sim |\xi_2|$ , and $\xi_1 , \xi_2 $ have different 
signs.\\
In this case we have $|\xi| = |\xi_1 - \xi_2| \sim 2|\xi_1| \sim 2|\xi_2|$ , 
and thus
\begin{eqnarray*}
C^p & \le & c \sup_{\xi,\sigma} \langle \xi \rangle^{lp} |\xi|^p \int d\xi_2 
d\sigma_2 \langle \sigma_1 \rangle^{-b_1 p} \langle \sigma_2 \rangle^{-b_1 p} 
|\xi|^{-2kp} \\
& \le &  c \sup_{\xi,\sigma} \langle \xi \rangle^{(l-2k)p} |\xi|^{p-1} \int  
\langle \sigma_1 \rangle^{-b_1 p} d\sigma_1 \int \langle \sigma_2 \rangle^{-b_1 
p} d\sigma_2 \, ,
\end{eqnarray*}
using $\frac{d\sigma_1}{d\xi_2} = 2\xi$ again. Under the assumption $l-2k+1 \le 
\frac{1}{p}$ $\Leftrightarrow$ $ 2k \ge l+\frac{1}{p'}$ this is finite. \\
{\bf Case 2:} $|\xi_1| \sim |\xi_2|$ , and $\xi_1,\xi_2$ have equal signs.\\
This implies $|\xi_1 + \xi_2| \sim 2|\xi_1| \sim 2|\xi_2|$ , and thus by 
(\ref{***}): 
$$ \xi(\xi_1+\xi_2) = (\xi_1 - \xi_2)(\xi_1 + \xi_2) = \sigma_1 - \sigma_2 - 
\sigma $$
and
$$ |\xi||\xi_2| \le c(|\sigma_1|+|\sigma_2|+|\sigma|) \, . $$
{\bf Case 2a:} $|\sigma_1|$ dominant ($\Rightarrow \, |\xi\xi_2| \le 
c|\sigma_1|$ ) ($|\sigma_2|$ dominant can be handled in the same way).\\
This also implies $|\xi| \langle \xi_2 \rangle \le c \langle \sigma_1 \rangle $ 
, which is evident for $|\xi_2| \ge 1$ , whereas $|\xi_2| \le 1$ implies $|\xi| 
= |\xi_1-\xi_2| \le c|\xi_2|$ (using $|\xi_1| \sim |\xi_2|$), so that $|\xi| 
\langle \xi_2 \rangle \le c \le c \langle \sigma_1 \rangle$. Thus
$$ C^p \le c \sup_{\xi,\sigma} \langle \xi \rangle^{lp} |\xi|^p \int d\xi_2 
d\sigma_2 \langle \xi_2 \rangle^{-2kp} \langle \sigma_1 \rangle^{-1-} \langle 
\sigma_2 \rangle^{-b_1 p} \langle \xi_2 \rangle^{1-b_1 p+} |\xi|^{1-b_1p +} \, 
. $$
Under the assumption $2k > \frac{1}{p} -b_1$ we have by $|\xi| \le c|\xi_2|$ 
and $\frac{d\sigma_1}{d\xi_2} = 2\xi$ :
$$ C^p \le c \sup_{\xi,\sigma} \langle \xi \rangle^{lp-2kp+1-b_1p+} 
|\xi|^{p-b_1p+} \int d\sigma_1 \langle \sigma_1 \rangle^{-1-} \int d\sigma_2 
\langle \sigma_2 \rangle^{-b_1p} \, . $$
Assuming $b_1 \le 1$ without loss of generality, this is finite, provided 
$lp-2kp+1-b_1p+p-b_1p<0$ $\Leftrightarrow$ $ l-2k+1 < 
-\frac{1}{p} +2b_1,$ which is fulfilled under our assumption $2k \ge l + 
\frac{1}{p'}$ $\Leftrightarrow$ $l-2k+1 \le \frac{1}{p}$ , because $b_1 > 
\frac{1}{p}$ .\\
{\bf Case 2b:} $|\sigma|$ dominant ( $\Rightarrow \, |\xi|\langle \xi_2 \rangle 
\le c \langle \sigma \rangle$ ). \\
We have
$$ C^p \le c \sup_{\xi,\sigma} \langle \xi \rangle^{lp} |\xi|^{p-cp} \int 
d\xi_2 d\sigma_2 \langle \xi_2 \rangle^{-2kp-cp} \langle \sigma_1 
\rangle^{-b_1p} \langle \sigma_2 \rangle^{-b_1p} \, . $$
Using $\langle \xi_2 \rangle^{-2kp-cp} \le c \langle \xi \rangle^{-2kp-cp}$ (by 
the assumption $2k \ge -c$) and $\frac{d\sigma_1}{d\xi_2} =2\xi$ , we get
$$ C^p \le c \sup_{\xi,\sigma} \langle \xi \rangle^{lp-2kp-cp} |\xi|^{p-cp-1} 
\int d\sigma_1 \langle \sigma_1 \rangle^{-b_1p} \int d\sigma_2 \langle \sigma_2 
\rangle^{-b_1p} \, . $$
The assumption $c \le \frac{1}{p'}$ implies $p-cp-1 \ge 0$ . Moreover we have 
$lp-2kp-cp+p-cp-1 \le 0$ $\Leftrightarrow$ $l-2k+1 \le \frac{1}{p} +2c$ , which 
is fulfilled under the assumption $l-2k+1 \le \frac{1}{p}$ $\Leftrightarrow$ $ 
2k \ge l + \frac{1}{p'}$ , so that $C^p$ is finite. \\
{\bf Case 3:} $|\xi_1| >> |\xi_2|$ ( $\Rightarrow \, |\xi| \sim |\xi_1|$ and 
$|\xi_2| << |\xi|$ ) (and similarly $|\xi_2| >> |\xi_1|$). \\
We have by $d\xi_2 = \frac{d\sigma_1}{2\xi}$ :
\begin{eqnarray*}
C^p & \le & c \sup_{\xi,\sigma} \int d\xi_2 d\sigma_2 \langle \xi 
\rangle^{(l-k)p} |\xi|^p \langle \xi_2 \rangle^{-kp} \langle \sigma_1 
\rangle^{-b_1p} \langle \sigma_2 \rangle^{-b_1p} \\
& \le & c \sup_{\xi,\sigma} \langle \xi \rangle^{(l-k)p} |\xi|^p \langle \xi 
\rangle^{-kp} |\xi|^{-1} \int d\sigma_1 \langle \sigma_1 \rangle^{-b_1p} \int 
d\sigma_2 \langle \sigma_2 \rangle^{-b_1p} \, ,
\end{eqnarray*}
which is finite, provided $l-2k+1 \le \frac{1}{p}$ $\Leftrightarrow$ $ 2k \ge 
l+\frac{1}{p'}$ .
\section{Proof of Theorem \ref{Theorem}}
{\bf Proof:} We construct a solution of the system of integral equations which 
belongs to our Cauchy problem by the contraction mapping principle. This can be 
achieved by using our Propositions \ref{Proposition 1.1}, \ref{Proposition 
1.2}, \ref{Proposition 1.3}, \ref{Proposition 1.4}, which give the necessary 
estimates for the nonlinearities, if one chooses $c_1 = 1-b_1-$ and $c=1-b-$. 
In this case the assumptions on the parameters in these propositions reduce to 
the assumptions in the theorem. We may apply Theorem \ref{Theorem 0.1} to our 
system, because its generalization from the case of a single equation to a 
system is evident.\\[1cm]

\end{document}